\documentclass[review]{elsarticle}

\usepackage{a4wide,amssymb,graphics,graphicx,textcomp}
\usepackage{enumerate,textcomp,multirow,amsmath}
\usepackage{algorithm}
\usepackage{algorithmic}
\usepackage{url}
\usepackage{xcolor} 
\usepackage{amsmath, amsthm}

\newtheorem{lemma}{Lemma}
\newtheorem{theorem}{Theorem}
\newtheorem{corollary}{Corollary}

\newtheorem{remark}{Remark}

\newcommand{\argmin}{\textrm{argmin}}
\newcommand {\mat}  [1] {\left[\begin{array}{#1}}
\newcommand {\rix}      {\end{array}\right]}

\newcommand{\eproof}{\space
    {\ \vbox{\hrule\hbox{\vrule height1.3ex\hskip0.8ex\vrule}\hrule}}\par}

\def\rank{\mathop{\mathrm{rank}}}

\newcommand{\C}{{\mathbb C}}
\newcommand{\R}{{\mathbb R}}

\definecolor{brightpink}{rgb}{1.0, 0.0, 0.5}

\journal{}

\begin{document}

\begin{frontmatter}

\title{On the non-symmetric semidefinite Procrustes problem}

\author{Mohit Kumar Baghel\fnref{iitd}}
\author{Nicolas Gillis\fnref{umons}}
\author{Punit Sharma\fnref{iitd}}
\fntext[iitd]{Department of Mathematics, Indian Institute of Technology Delhi, Hauz Khas, 110016, India; \texttt{\{maz188260, punit.sharma\}@maths.iitd.ac.in.}
MB acknowledges the support of institute Ph.D. fellowship by IIT Delhi, India.  PS acknowledges the support of the DST-Inspire Faculty Award (MI01807-G) by Government of India and FIRP project (FIRP/Proposal Id - 135) by IIT Delhi, India.
}
\fntext[umons]{Department of Mathematics and Operational Research,
Facult\'e Polytechnique, Universit\'e de Mons, Rue de Houdain~9, 7000 Mons, Belgium; \texttt{nicolas.gillis@umons.ac.be.} 
NG acknowledges the support by ERC starting grant No 679515, the Fonds de la Recherche Scientifique - FNRS
and  the  Fonds  Wetenschappelijk  Onderzoek - Vlaanderen (FWO) 
under EOS project O005318F-RG47. 
}

%
%

\begin{abstract}
In this paper, we consider the non-symmetric positive semidefinite Procrustes (NSPSDP) problem: Given two matrices $X,B \in \mathbb{R}^{n,m}$, find the matrix $A \in \mathbb{R}^{n,n}$ that minimizes the Frobenius norm of $AX-B$ and which is such that $A+A^T$ is positive semidefinite. We generalize the semi-analytical approach for the symmetric positive semidefinite Procrustes problem, where $A$ is required to be positive semidefinite, that was proposed by Gillis and Sharma (A semi-analytical approach for the positive semidefinite Procrustes problem, Linear Algebra Appl. 540, 112-137, 2018). As for the symmetric case, we first show that the NSPSDP problem can be reduced to a smaller NSPSDP problem that always has a unique solution and where the matrix $X$ is diagonal and has full rank. Then, an efficient semi-analytical algorithm to solve the NSPSDP problem is proposed, solving the smaller and well-posed problem with a fast gradient method which guarantees a linear rate of convergence. This algorithm is also applicable to solve the complex NSPSDP problem, where $X,B \in \mathbb{C}^{n,m}$,  as we show the complex NSPSDP problem can be written as an overparametrized real NSPSDP problem. The efficiency of the proposed algorithm is illustrated on several numerical examples. 
\end{abstract}


\end{frontmatter}


\section{Introduction}

For given $X$, $B \in \C^{n,m}$, and $\mathcal S \subseteq \C^{n,n}$, the mapping problem associated with the set $\mathcal S$ is to find $
\Delta \in \mathcal S$ such that $\Delta X=B$, see~\cite{MacMT08,Adh08} for various structured mapping problems.
 When a mapping problem is inconsistent (i.e., there does not exist $\Delta \in \mathcal S$ such that $\Delta X=B$), it is natural to look for $\Delta $ such that 
$\Delta X$ is as close to $B$ as possible, and solve 
\begin{equation}\label{eq:leastsquare}
\min_{\Delta \in \mathcal S} {\| \Delta X-B\|},
\end{equation}
where $\|\cdot\|$ is some matrix norm. 
Problems of the form~\eqref{eq:leastsquare} have been widely studied and are known as Procrustes problems (they can also be seen as least square problems). 
For a given structure, the first aim to solve the problem~\eqref{eq:leastsquare} analytically (if possible), if not then the aim is to find a computable algorithmic solution with the desired accuracy.
By considering different structures on $\mathcal S$, a variety of Procrustes problems
can be defined, see~\cite{Gre52,Hig88c,Sch66,Woo96,JinWPC19,KrisLVP04,Kris03,GilS18}. 

Such Procrustes problems occur in the estimation of compliance or stiffness  matrix in solid and elastic structures that are based on observations of force and displacement. The constraint on the solution of the problem depends on the physics of the problem. For example the stiffness matrix is symmetric in determination of space craft altitudes~\cite{Hig88c,Lar66}, positive semidefinite in elastic structure~\cite{SufH93,Bro68}, and orthogonal in computer graphics, multidimensional scaling and factor analysis~\cite{GowD04,Gre52}. Unlike the stiffness matrix in elastic structure, the measured local stiffness matrices in deformable models are non-symmetric positive semidefinite~\cite{KrisLVP04}. 
A matrix $A\in R^{n,n}$ is called non-symmetric positive semidefinite (NSPSD) if $A+A^T$ is positive semidefinite.

 In this paper, we consider the following  non-symmetric positive semidefinite Procrustes (NSPSDP) problem:

\textbf{ NSPSDP problem}:~
Given $n \times m$ real matrices $X$ and $B$, the NSPSDP-problem is defined by
\begin{equation} \label{nspsdproc}
\inf_{ A \in \mathcal{N}_{\succeq}^{S_n}  } {\|AX-B\|}_F^2, \tag{$\mathcal{P}$}
\end{equation}
where $\mathcal{N}_{\succeq}^{S_n} $ is the set of all NSPSD matrices of size $n \times n$, and
${\|\cdot\|}_F$ stands for the Frobenius norm. 
This problem was first defined in~\cite{KrisLVP04}. It is also known as the dissipative mapping problem in the exact case, that is, when one is looking for $A$ such that $AX = B$; 
see~\cite{BhaGS21} and the references therein. 

 The solutions to the NSPSDP problem can be useful in the estimation of local compliance matrices during deformable object  in various engineering applications in robotics, computer graphics and medical simulation~\cite{KrisLVP04}. The deformable object to be  modelled is  a passive object, i.e., it does not generate energy in deformation. Passivity of a deformable object restricts its compliance at any surface point such that the traction generated from an applied displacement is opposed to the applied traction. This requirement can be expressed as $p^TE p \geq 0$, where $E$ is the Green's functions matrix and $p$ is the point load. In estimation of $E$ under noisy measurements of traction and displacements, one has to ensure that the 
 compliance matrix $E$ satisfies this constraint, i.e., $E$ is NSPSD.

A closely related problem to~\eqref{nspsdproc} is the positive semidefinte procrustes (PSDP) problem, i.e.,   when $\mathcal S$ in~\eqref{eq:leastsquare} is the set of positive semidefinte (PSD) matrices. The PSDP problem has been extensively studied in~\cite{Bro68,Woo96,Kris03,GilS18,JinWPC19}. Recently, a semi-analytical approach was proposed in~\cite{GilS18} to solve the PSDP problem where the original problem was reformulated to an equivalent PSDP problem with diagonal $X$, and then a fast gradient method was proposed to solve this reformulated problem algorithmically. An inner product based approach to solve PSDP problem is attempted in~\cite{JinWPC19} which solved the problem by computing the exact solution analytically.  
In comparison with the PSDP problem, the NSPSDP problem has attracted much less attention~\cite{KrisLVP04,Kris03}. In~\cite{KrisLVP04}, authors proposed an algorithm based on an interior-point method for solving the NSPSDP problem and applied this method to the estimation of local compliance matrices of elastic solids. 

 Motivated by~\cite{GilS18}, we consider the NSPSDP problem~\eqref{nspsdproc} and  present a semi-analytical approach to solve  NSPSDP problem with solutions of desired accuracy. As far as we know, this is the first attempt in solving the NSPSDP problem analytically.  
 
 A matrix $A \in \C^{n,n}$ is called non-Hermitian positive semidefinite (NHPSD) if $A+A^*$ is positive semidefinite. In this paper, we also consider the complex version of the problem~\eqref{nspsdproc}. More precisely, 
for given $n \times m$ complex matrices $X$ and $B$, we aim to solve the non-Hermitian positive semidefinite procrustes (NHPSDP) problem  
\begin{equation} \label{complex_main_problem} 
\inf_{A  \in \mathcal{N}_\succeq^{\mathcal{H}_n}} {\| AX-B\|}_F^2  , \tag{$\mathcal{P_\C}$}
\end{equation}
where $\mathcal{N}_\succeq^{\mathcal{H}_n}$ is the set of all NHPSD matrices of size $n \times n$.


\subsection{Contributions and outline of the paper}

In Section~\ref{sec:prelim}, we present some preliminary results that will be needed to solve the NSPSDP problem.  
In Section~\ref{sec:nspsdp}, we consider the NSPSDP problem~\eqref{nspsdproc}. A semi-analytic approach was proposed in~\cite{GilS18} for solving PSDP problem. Motivated by~\cite{GilS18}, we present a semi-analytic approach to solve NSPSDP problem where we reduce the original problem~\eqref{nspsdproc} into a smaller NSPSDP problem with a full rank diagonal $X$ that always has a unique solution. 
Using the solution to this smaller problem, we present  a family of NSPSD matrices that either solve the problem~\eqref{nspsdproc} when the infimum is attained or give an approximate solution arbitrary close to the infimum when it is not attained. We then employ an optimal first-order method to solve the smaller problem algorithmically in Section~\ref{sec:algo}. In Section~\ref{sec:compnspsdp}, we consider the NHPSDP problem~\ref{complex_main_problem} and transform this problem into a real 
NSPSDP problem of double size, for which the algorithm proposed in Section~\ref{sec:algo} is applicable. 
 In Section~\ref{sec:numerical}, we illustrate the performance of our proposed algorithm on some numerical examples and compare the results with the algorithm presented in~\cite{KrisLVP04}. 
 
\paragraph{Notation} In the following, we denote the identity matrix of size $n \times n$ by $I_n$, the spectral norm of a matrix or a vector by $\|\cdot\|$ and the Frobenius norm by ${\|\cdot\|}_F$. The Moore-Penrose pseudoinverse of a matrix or a vector $X$ is denoted by $X^\dagger$
and  $\mathcal P_X=I_n-XX^\dagger$ denotes the orthogonal projection onto the null space of $n \times n$ matrix $X^*$. For a square matrix $A$, its Hermitian and skew-Hermitian parts are respectively denoted by $A_H=\frac{A+A^*}{2}$ and $A_S=\frac{A-A^*}{2}$.
For $A=A^* \in \mathbb F^{n,n}$, where $\mathbb F \in \{\R,\C\}$, we denote
$A \succ 0$ ($A\prec 0$) and $A\succeq 0$ ($A\preceq$) if $A $ is Hermitian positive definite (negative definite) and Hermitian positive semidefinite (negative semidefinite).  $\Lambda(A)$ denotes the set of all eigenvalues of the matrix $A$. We follow the notation from~\cite{HorJ94}
 and denote any complex matrix $A \in \C^{n,n}$ by $A=A_1+iA_2$, where 
 $A_1,A_2 \in \R^{n,n}$, and define $R(A):=\mat{cc}A_1 &A_2 \\-A_2 & A_1\rix$.

\section{Preliminary results} \label{sec:prelim}

In this section, we recall some preliminary results from the literature
and also derive some basic results in the form which will be useful in the later sections. 

The following lemma which gives equivalent criteria to check for a block matrix to be positive semidefinite, will be used to solve the  NSPSDP problem.

\begin{lemma}{\rm \cite{Alb69}} \label{lem:psd_character}
Let the integer $s$ be such that $0<s<n$, and $R=R^T \in \R^{n,n}$ be partitioned as
$R=\mat{cc}B & C^T\\C & D\rix$ with $B\in \R^{s,s}$, $C\in \R^{n-s,s}$ and $D \in \R^{n-s,n-s}$. Then $R \succeq 0$ if and only if
\begin{enumerate}
\item $B \succeq 0$,
\item $\operatorname{null}(B) \subseteq \operatorname{null}(C)$, and
\item $D-CB^{\dagger}C^T \succeq 0$, where $B^{\dagger}$ denotes the Moore-Penrose pseudoinverse of $B$.
\end{enumerate}
\end{lemma}
The following lemma is analogous to~\cite[Lemma 4]{GilS18} for NSPSD 
matrices which will be used
in finding minimal Frobenius norm solution to the NSPSDP problem.
\begin{lemma}\label{lem:minnorm}
Let the integer $s$ be such that 
$0<s<n$. Let $B \in \mathcal{N}_{\succeq}^{S_s}$ and $C \in \R^{n-s,s}$ be such that $\text{{\rm null}}(B_H)\subseteq \text{{\rm null}}(\frac{C}{2})$. Define 
\[
\mathcal M:=\left\{(K,R)\in (\R^{n-s,n-s})^2:~K-\frac{1}{4}CB_H^\dagger C^T \succeq 0,\,R^T=-R\right \},
\]
 and  $g:\mathcal M \longmapsto \R^{n,.n}$ by
 \[
 g(K,R)= \mat{cc}B&0\\C& K+R\rix.
 \]
 Then the matrix pair $(\widetilde K, \widetilde R):=(\frac{1}{4}CB_H^\dagger C^T,0) \in \mathcal M$ is the unique solution of the minimal Frobenius norm problem, that is,
 \[ 
\min_{(K,R) \in \mathcal M}{\|g(K,R)\|}_F={\|g(\widetilde K, \widetilde R)\|}_F.
\]
\end{lemma}
\proof The proof is similar to~\cite[Lemma 4]{GilS18} and hence left for the reader.
\eproof
The trace of product of a symmetric positive semidefinite and a non-symmetric positive semidefinite matrix is always nonnegative, as shown in the following lemma. 
\begin{lemma}\label{lem:trace}
Let $A \in \R^{n,n}$ such that $A\succeq 0$ and let $B \in \mathcal{N}_{\succeq}^{S_n}$. Then $\text{{\rm trace}}(AB) \geq 0$.
\end{lemma}
\proof 
First note that since $A \succeq 0$, we have $A=UDU^T$ for some orthogonal
matrix $U \in \R^{n,n}$ and a diagonal matrix $D \in \R^{n,n}$ with all nonnegative diagonal entries. Thus we have
\begin{equation}\label{eq:lemtr1}
\text{trace}(AB_S)=\text{trace}(UDU^TB_S)=\text{trace}(DU^TB_SU)=0,
\end{equation}
since $U^TB_SU$ is skew-symmetric with diagonal entries all equal to zero. This implies that 
\begin{eqnarray}\label{eq:lemtr2}
\text{trace}(AB)=\text{trace}(A(B_H+B_S))=\text{trace}(AB_H)+\text{trace}(AB_S)=\text{trace}(AB_H)\geq 0,
\end{eqnarray}
where the last equality follows from~\eqref{eq:lemtr1} and due to  the fact that 
$A \succeq 0$ and $B_H \succeq 0$, imply that all eigenvalues of the product 
$AB_H$ are nonnegative~\cite[Lemma 3]{GilS18}.
\eproof

The following result guarantees a unique solution to NSPSDP problem when $X$ is of full column rank. 
\begin{theorem}{\rm \cite[Theorem 2.4.6]{Kris03}}\label{lem:KrisNSPSDP}
 If $X$ has full column rank, then the NSPSDP problem~\eqref{nspsdproc} has a unique optimal solution.
\end{theorem}
The next lemma gives an equivalent real condition for a complex matrix  
$A \in \mathcal{N}_\succeq^{\mathcal{H}_n}$, which is a generalization of 
~\cite[Lemma-3.1]{KisH07} for NHPSD matrices.

\begin{lemma} \label{complex_real_nspsd_lemma}
	Let $A=A_r + i A_j \in \mathbb{C}^{n,n}$ with $A_r, A_j \in \mathbb{R}^{n,n}$. Then  $A \in \mathcal{N}_\succeq^{\mathcal{H}_n}$ if and only if  $P:= \mat{cc}
		A_r & A_j\\
		-A_j & A_r
		\rix \in \mathcal{N}_\succeq^{\mathcal{S}_{2n}}$. 
\end{lemma}
\proof The proof follows by the following observation: for any $x=x_r+ix_j \in \C^n$, where $x_r,x_j \in \R^n$, we have
\begin{eqnarray} \label{logic1}
x^*(A + A^*)x&=& x^*((A_r + i A_j) + (A_r + i A_j)^*) x     \nonumber \\  
&=& (x_r + i x_j)^* (A_r + A_r ^T) (x_r + i x_j) + (x_r + i x_j)^* (A_j - A_j ^T) (x_r + i x_j) \nonumber \\
&=& (x_r ^T - i x_j ^T) (A_r + A_r ^T) (x_r + i x_j)+ i(x_r ^T - i x_j ^T) (A_j - A_j ^T) (x_r + i x_j) \nonumber \\
&=& x_r^T (A_r + A_r ^T) x_r +  x_j ^T (A_r + A_r ^T) x_j + i x_r^T (A_r + A_r ^T) x_j - i x_j ^T (A_r + A_r ^T) x_r  \nonumber \\
&&+ i x_r ^T (A_j - A_j ^T) x_r + i x_j^T (A_j - A_j ^T) x_j - x_r ^T (A_j - A_j ^T) x_j + x_j ^T (A_j - A_j ^T) x_r  \nonumber \\
&=& x_r^T (A_r + A_r ^T) x_r  + x_j ^T (A_r + A_r ^T) x_j
-x_r^T (A_j - A_j ^T) x_j+ x_j ^T (A_j - A_j ^T) x_r ,
\end{eqnarray}
where  the last identity holds because $x_r^T (A_r + A_r ^T) x_j =x_j^T (A_r + A_r ^T) x_r $ since $A_r + A_r ^T$ is symmetric,  and $x_r ^T (A_j - A_j ^T) x_r =0$ and $x_j^T (A_j - A_j ^T) x_j= 0$ since $A_j-A_j^T$ is skew-symmetric. Also, 
for $\widetilde x= [-x_r^T ~\,x_j^T]^T $, we have 
\begin{eqnarray} \label{logic2}
\widetilde{x}^T (P + P^T)\widetilde{x}
&=&\begin{bmatrix} -x_r^T & x_j^T\end{bmatrix}
\begin{bmatrix}
A_r + A_r^T & A_j -A_J^T\\-A_j + A_j^T & A_r + A_r^T\end{bmatrix} 
\begin{bmatrix}-x_r \\x_j \end{bmatrix} \nonumber \\
&=& x_r^T(A_r + A_r^T)x_r + x_j^T(A_j - A_j^T)x_r- x_r^T(A_j-A_j^T)x_j + x_j^T(A_r + A_r^T) x_j. \nonumber \\
\end{eqnarray}
Clearly, from~\eqref{logic1} and~\eqref{logic2}, $x^*(A + A^*)x \geq 0$ if and only if $\widetilde{x}^T (P + P^T)\widetilde{x} \geq 0$.
\eproof

The following lemma gives the projection $P_{\succeq}(Z)$ of $Z \in \C^{n,n}$ onto the cone of complex positive semidefinite matrices. 

\begin{lemma}{\rm \cite{Hig88c}}\label{lem:com5}
Let $A \in \C^{n,n}$. Then 
\[
 P_{\succeq} (A):={\rm argmin}_{X \in \C^{n,n},\, X \succeq 0} {\|A-X\|}_F^2 = V\left(\max{(\Gamma,0)}\right)V^*,
\]
where $V \Gamma V^*$ is an eigenvalue decomposition of the Hermitian matrix $\frac{A+A^*}{2}$.
\end{lemma}

We close this section with two elementary lemmas which will be used in Section~\ref{sec:compnspsdp} for complex NHPSDP problem. 




\begin{lemma}{{\rm \cite[Problems 1.3.P20]{HorJ94}}}\label{lem:com1} 
Let $A=A_1+iA_2$, where   $A_1,A_2 \in \R^{n,n}$ and let $U=\frac{1}{\sqrt{2}}\mat{cc}I_n &iI_n \\iI_n & I_n\rix$. Then $U^*=\overline U=U^{-1}$, and 
$U^*R(A)U=\mat{cc}A &0 \\ 0 & \overline A \rix$, where $\overline A$ is the complex conjugate of $A$ and $R(A)=\mat{cc}A_1 & A_2\\-A_2 & A_1\rix$. 
\end{lemma}

\begin{lemma}\label{lem:com2}
Let $X=X_1+iX_2 \in \C^{n,n}$, where $X_1,X_2 \in \R^{n,n}$. Then {{$X \succeq 0$}} if and only if $R(X) \succeq 0$. Similarly, $X \in \mathcal{N}_\succeq ^{\mathcal H_{n}}$ if and only if $R(X) \in \mathcal{N}_\succeq ^{S_{2n}}$.
\end{lemma}

\section{ NSPSDP problem: A semi-analytical solution } \label{NSPSD Problem} \label{sec:nspsdp}
In this section, we present a semi-analytical method to solve the NSPSDP problem. We first reduce the original problem~\eqref{nspsdproc} into a smaller size problem with full rank diagonal $X$ that always has a unique solution. 
Then, assuming the solution to this smaller problem is known, we provide a characterization for all the solutions to the NSPSDP problem~\eqref{nspsdproc} when the infimum is  attained, or provide a family of solutions that achieve the desired accuracy when the infimum is not attained. 
This reduction is achieved in the following result which is analogous to~\cite[Theorem 1]{GilS18}. 
An optimal first-order method is proposed 
to solve the smaller problem algorithmically in Section~\ref{sec:algo}.

 	\begin{theorem} \label{thetheorem}
		Let $X,B \in \mathbb{R}^{n,m},$ and let $r=\rank(X).$ Suppose that $X=U \Sigma V^T$ is a singular value decomposition of $X,$ where $U=[U_1 ~ U_2] \in \mathbb{R}^{n,n}$ with $U_1 \in \mathbb{R}^{n,r}, V=[V_1 ~ V_2] \in \mathbb{R}^{m,m}$ with $V_1 \in \mathbb{R}^{m,r},$ and $\Sigma= \begin{bmatrix}
\Sigma_1 & 0\\0 & 0 \end{bmatrix} \in \mathbb{R}^{n,m}$ with $\Sigma_1 \in \mathbb{R}^{r,r}$. Then 
\begin{equation} \label{semi_analytic_eq}
\inf_{A \in \mathcal{N}_{\succeq}^{S_n}} {\|AX-B\|}^2 = \min_{A_{11} \in \mathcal{N}_{\succeq}^{S_r}} {\|A_{11}\Sigma_1 - U_1^T B V_1\|}_F ^2 + {\|B V_2\|}_F ^2.
\end{equation} 
Further, let $\hat{A}_{11} \in \mathcal{N}_{\succeq}^{S_r}$ be such that 
\begin{equation}\label{assumtion}
\hat{A}_{11} = \text{{\rm argmin}}_{A_{11} \in  \mathcal{N}_{\succeq}^{S_r}} {\|A_{11} \Sigma_1 - U_1^T B V_1\|}_F ^2,
\end{equation}
and let $\hat{H}_{11}=\frac{1}{2}(\hat{A}_{11}+\hat{A}_{11}^T)$ and $\hat{S}_{11}=\frac{1}{2}(\hat{A}_{11}-\hat{A}_{11}^T)$.
Then the following hold.
\begin{enumerate}
\item[(i)]~If $\text{{\rm null}}(\hat{H}_{11}) \subseteq \text{{\rm null}}(\frac{1}{2}U_2^{T} B V_1 \Sigma_1 ^{-1}) $, then $A_{opt}$ attains the infimum in~\eqref{semi_analytic_eq} if and only if 
\begin{equation} \label{A_opt}
A_{opt} := U_1 \hat{A}_{11} U_1 ^{T} + U_2(U_2 ^T B V_1 \Sigma_1 ^{-1}) U_1 ^{T} + U_2 K U_2^{T}+ U_2RU_2^T,
\end{equation}
where $K \in \mathbb{R}^{n-r,n-r}$ is such that $K-\frac{1}{4}(U_2^{T} B V_1 \Sigma_1 ^{-1}) (\hat H_{11}) ^{\dagger} (U_2 ^T B V_1 \Sigma_1 ^{-1})^T \succeq 0$ and $R \in \R^{n-r,n-r}$ such that  $R^T=-R$.		

\item[(ii)]~Otherwise, the infimum in \eqref{semi_analytic_eq} is not attained. 
Let $\rank(\hat{H}_{11})=s < r$ and let $\epsilon > 0$ be sufficiently small according to \eqref{epsilon_inequality_eq}. Let 
${\hat H}_{11}= \begin{bmatrix}\hat{U}_1 & \hat{U}_2\end{bmatrix} 
\begin{bmatrix}\hat{\Sigma}_1 & 0\\0 & 0\end{bmatrix} 
\begin{bmatrix}\hat{U}_1^T\\\hat{U}_2^{T} \end{bmatrix}$ 
be a SVD of $\hat{H}_{11},$ where $\hat{U}_1 \in \mathbb{R}^{r,s}$ and $\hat{\Sigma}_1 \in \mathbb{R}^{s,s}.$ Define 
\begin{equation}\label{A11_epsilon} 
\hat{H}_{11}^{\epsilon} :=  \begin{bmatrix} \hat{U}_1 & \hat{U}_2\end{bmatrix} 
\begin{bmatrix}\hat{\Sigma}_1 & 0\\0 & \Gamma\end{bmatrix} 
\begin{bmatrix} \hat{U}_1^T\\ \hat{U}_2^{T}\end{bmatrix},
\end{equation}
where $\Gamma \in \mathbb{R}^{r-s,r-s}$ is a diagonal matrix with diagonal entries each equal to $\frac{\epsilon}{\beta}$, where 
\begin{equation}\label{eq:def_beta}
\beta = \begin{cases} 
4 \sqrt{(r-s)} {\|\Sigma_1\|}_F {\|\hat{A}_{11} \Sigma_1 - U_1^T B V_1\|}_F & if {\|\hat{A}_{11} \Sigma_1 - U_1^T B V_1\|}_F \neq 0,  \\
4 \sqrt{(r-s)} {\|\Sigma_1\|}_F  & otherwise.
			\end{cases} 
\end{equation}
Define ${\hat A}_{11}^\epsilon:={\hat H}_{11}^\epsilon+{\hat S}_{11}$ and 
\begin{equation} \label{A_eps}
A_{\epsilon} := U_1 \hat{A}_{11}^{\epsilon} U_1^{T} + U_2(U_2^T B V_1 \Sigma_1^{-1}) U_1^{T} +  U_2 K_{\epsilon} U_2^T+U_2 R U_2^T,
\end{equation}
where $K_\epsilon \in \mathbb{R}^{n-r,n-r}$ is such that $K_{\epsilon} -\frac{1}{4}
(U_2^{T} B V_1 \Sigma_1^{-1}) (\hat{H}_{11} ^{\epsilon})^{-1} (U_2 ^T B V_1 \Sigma_1 ^{-1})^T \succeq 0$, and $R \in \R^{n-r,n-r}$ such that  $R^T=-R$. 
		Then $A_{\epsilon} \in \mathcal{N}_{\succeq}^{S_n}$ and 
		\begin{equation}
			{\| A_{\epsilon} X - B\|}_F ^{2} < \inf_{A \in \mathcal{N}_{\succeq}^{S_n} } {\| A X -B\|}_F ^{2} + \epsilon.
		\end{equation}
		\end{enumerate}		
\end{theorem}
\proof Let $A \in \mathbb{R}^{n,n}$ and set 
\begin{equation} \label{svd_eq}
\hat{A}:= U^T A U= \hat A_H + \hat A_S,
\end{equation}
where 
\[
 \hat A_H=U^T\left(\frac{A+A^T}{2}\right)U
=  \begin{bmatrix}
H_{11} & H_{21}^T \\
H_{21} & H_{22}
\end{bmatrix},\quad \text{and}\quad 
 \hat A_S=U^T\left(\frac{A-A^T}{2}\right)U
=
\begin{bmatrix}
S_{11} & -S_{21}^T \\
S_{21} & S_{22}
\end{bmatrix},
\]
where $H_{11},\,S_{11} \in \mathbb{R}^{r,r}$, 
$ H_{21},S_{21} \in \mathbb{R}^{n-r,r}$,  $H_{22},\,S_{22} \in \mathbb{R}^{n-r,n-r}$. Since $A+A^T \succeq 0$ if and only if $\hat A_H \succeq 0$. This implies by Lemma~\ref{lem:psd_character}  that  $\hat A_H \succeq 0$ if and only if
$H_{11}\succeq 0$, $\text{null}(H_{11})\subseteq \text{null}(H_{21})$ and $H_{22}-H_{21}H_{11}^\dagger H_{21}^T \succeq 0$. Also $A-A^T$ is skew-symmetric if and only if $S_{11}^T=-S_{11}$ and $S_{22}^T=-S_{22}$. Thus we have 
\begin{eqnarray}\label{eq:proof1}
{\|AX-B\|}_F^2   &=& { \| U \hat A U^T X - B\|}_F ^2
= {\|\hat A U^TX-U^TB\|}_F^2\nonumber \\
&=&  {\|(\hat A_H+\hat A_S) U^TX-U^TB\|}_F^2
= {\left\| (\hat A_H+\hat A_S) \mat{c}U_1^TX \\ 0 \rix -\mat{c}U_1^TB \\U_2^TB\rix\right\|}_F^2  \nonumber\\
&=& {\left\| \mat{c}(H_{11}+S_{11})U_1^TX-U_1^TB \\ (H_{21}+S_{21})U_1^TX-U_2^TB  \rix \right\|}_F^2 \nonumber \\
&=& {\left\|(H_{11}+S_{11})U_1^TX-U_1^TB\right\|}_F^2 +  {\left\|(H_{21}+S_{21})U_1^TX-U_2^TB \right\|}_F^2 \nonumber \\
&=& {\left\|(H_{11}+S_{11})U_1^T\mat{cc}U_1 &U_2 \rix \mat{cc}\Sigma_1 &0\\0&0\rix-U_1^TB \mat{cc} V_1 & V_2 \rix\right\|}_F^2 \nonumber \\
&+& {\left\|(H_{21}+S_{21})U_1^T\mat{cc}U_1 &U_2 \rix \mat{cc}\Sigma_1 &0\\0&0\rix-U_2^TB \mat{cc} V_1 & V_2 \rix\right\|}_F^2 \nonumber\\
&=& {\left\|(H_{11}+S_{11})\mat{cc}I &0 \rix \mat{cc}\Sigma_1 &0\\0&0\rix- \mat{cc} U_1^TBV_1 & U_1^TBV_2 \rix\right\|}_F^2 + \nonumber\\
&& {\left\|(H_{21}+S_{21})\mat{cc}I &0 \rix \mat{cc}\Sigma_1 &0\\0&0\rix- \mat{cc} U_2^TB V_1 & U_2^TB V_2 \rix\right\|}_F^2 \nonumber \\
&=& {\left \| \mat{cc}(H_{11}+S_{11}) \Sigma_1 -U_1^TBV_1 & -U_1^TBV_2\rix\right \| }_F^2 + \nonumber\\
&&{\left \| \mat{cc}(H_{21}+S_{21}) \Sigma_1 -U_2^TBV_1 & -U_2^TBV_2\rix\right \| }_F^2  \nonumber \\
&=& {\left \| (H_{11}+S_{11}) \Sigma_1 -U_1^TBV_1\right \|}_F^2 + {\left \| U_1^TBV_2\right \|}_F^2 +\nonumber \\
&& {\left \|(H_{21}+S_{21}) \Sigma_1 -U_2^TBV_1 \right \|}_F^2 + 
{\left \|U_2^TBV_2 \right \|}_F^2 \nonumber \\ 
&=& {\left \| (H_{11}+S_{11}) \Sigma_1 -U_1^TBV_1\right \|}_F^2  + {\left \|(H_{21}+S_{21}) \Sigma_1 -U_2^TBV_1 \right \|}_F^2 + 
{\left \|BV_2 \right \|}_F^2, \nonumber \\
\end{eqnarray}
where the last equality follows by the fact that 
\[
\left\|U_1^TBV_2\right\|_F^2 + \left\|U_2^TBV_2\right\|_F^2 
= \left\| \left[U_1^T \; U_2^T \right]   BV_2\right\|_F^2 
= \left\|BV_2\right\|_F^2 
\] 
since ${\|\cdot\|}_F$ is unitary invariant. By taking infimum in~\eqref{eq:proof1} over $\mathcal{N}_{\succeq}^{S_n}$, we obtain
\begin{eqnarray} \label{inf_eq}
\inf_{A \in \mathcal{N}^{S_n}_{\succeq}} {\|AX-B\|}_F^2 
=\inf_{H_{11},S_{11} \in \R^{r,r}, H_{21},S_{21} \in \R^{n-r,r}}   \Big \{ {\|(H_{11}+S_{11}) \Sigma_1 - U_1^T B V_1\|}_F^2  \nonumber \\+ 
{\|(H_{21}+S_{21}) \Sigma_1 - U_2^T B V_1\|}_F^2 + {\|BV_2\|}_F^2 \Big \}
\end{eqnarray}
such that $H_{11} \succeq 0, \text{null}(H_{11}) \subseteq \text{null}(H_{21}),S_{11}^T=-S_{11}$. Note that in~\eqref{inf_eq},  the infimum is independent of the choice of 
$H_{22}$ and $S_{22}$, and therefore the conditions 
 $H_{22}-H_{21}H_{11}^\dagger H_{21}^T \succeq 0$  and $S_{22}^T=-S_{22}$
 can be removed. In view of~\eqref{inf_eq}, we have 
 \begin{eqnarray} \label{inequality_eq1}
 \inf_{A \in \mathcal{N}^{S_n}_{\succeq}} {\|AX-B\|}_F^2 \geq 
 \inf_{H_{11} \succeq 0,S_{11}^{T}=-S_{11}} {\|(H_{11}+S_{11}) \Sigma_1 - U_1^T B V_1\|}_F^2 \nonumber \\
+ \inf_{K \in \mathbb{R}^{n-r,r}}{\|K \Sigma_1 - U_2^T B V_1\|}_F^2 + {\|BV_2\|}_F^2,
 \end{eqnarray}
since dropping the constraint $\text{null}(H_{11}) \subseteq \text{null}(H_{21}) $
in~\eqref{inf_eq} will only result in a smaller value of the infimum. This implies that 
\begin{eqnarray}\
 \inf_{A \in \mathcal{N}^{S_n}_{\succeq}} {\|AX-B\|}_F^2 &\geq&  
  \inf_{A_{11} \in \mathcal{N}^{S_r}_{\succeq}} {\|A_{11} \Sigma_1 - U_1^T B V_1\|}_F^2+ {\|BV_2\|}_F^2 \label{eq:proof2}\\
  &=& {\|\hat A_{11} \Sigma_1 - U_1^T B V_1\|}_F^2  + {\|BV_2\|}_F^2, \nonumber
\end{eqnarray}
 where the last equality holds since (i) the value of the second infimum in~\eqref{inequality_eq1} is zero as it is attained at $K=U_2^TBV_1\Sigma_1^{-1}$, and (ii) by Lemma~\ref{lem:KrisNSPSDP}, the first infimum is attained at a unique $\hat A_{11} \in \mathcal{N}^{S_r}_{\succeq}$. Define $\hat H_{11}=\frac{1}{2}(\hat A_{11}+\hat A_{11}^T)$ and $\hat S_{11}=\frac{1}{2}(\hat A_{11}-\hat A_{11}^T)$. As done in~\cite{GilS18}, we consider two cases to show the equality in~\eqref{eq:proof2}.
		
\underline{Case-1:} $\text{null}(\hat{H}_{11}) \subseteq \text{null}(\frac{1}{2}U_2^T BV_1 \Sigma_1^{-1})$.
 Let $H_{21}=\frac{1}{2}U_2^T BV_1 \Sigma_1^{-1}$, $S_{21}=\frac{1}{2}U_2^T BV_1 \Sigma_1^{-1}$,
 $H_{11}=\hat{H}_{11}$ and $S_{11}=\hat S_{11}$ in \eqref{svd_eq}, and define
\begin{equation*}
A_{opt} := U \left( \begin{bmatrix}\hat{H}_{11} & \frac{1}{2}(U_2^T BV_1 \Sigma_1^{-1})^T\\ \frac{1}{2}U_2^T BV_1 \Sigma_1^{-1} & K \end{bmatrix}
+ \begin{bmatrix} \hat{S}_{11} & -\frac{1}{2}(U_2^T BV_1 \Sigma_1^{-1})^T\\ \frac{1}{2}U_2^T BV_1 \Sigma_1^{-1} & R \end{bmatrix} \right)U^T , 
\end{equation*}
for some $K \in \mathbb{R}^{n-r,n-r}$ such that 
$K-\frac{1}{4}(U_2^{T} B V_1 \Sigma_1 ^{-1}) \hat{H}_{11} ^{\dagger} (U_2 ^T B V_1 \Sigma_1 ^{-1})^T \succeq 0$ and $R \in \R^{n-r,n-r}$ such that $R^T=-R$, which upon simplification gives \eqref{A_opt}. 
Thus by Lemma~\ref{lem:psd_character}, $A_{opt} \in \mathcal{N}_{\succeq}^{S_n}$ and from~\eqref{eq:proof1}, it satisfies 
 \begin{equation}
 {\|A_{opt}X-B\|}_F^2 ={\|(\hat{H}_{11} + \hat{S}_{11}) \Sigma_1- U_1^T B V_1\|}_F^2 + {\|BV_2\|}_F^2.
 \end{equation}
This implies equality in~\eqref{eq:proof2} and hence proves~(i).

 \underline{Case-2:} $\text{null}(\hat{H}_{11}) \not \subseteq \text{null}(\frac{1}{2}U_2^T BV_1 \Sigma_1^{-1})$. Let $s=\text{rank}(\hat H_{11})$ and $\epsilon$
 sufficiently small such that
\begin{equation} \label{epsilon_inequality_eq}
 0 < \epsilon < \begin{cases}
 \min\{1,\,{\|\hat{A}_{11} \Sigma_1- U_1^T B V_1\|}_F^2\} & \text{if} ~{\|\hat{A}_{11} \Sigma_1- U_1^T B V_1\|}_F \neq 0,\\
 1 & \text{otherwise}.
 \end{cases}
 \end{equation}
Let  $\hat H_{11}=\mat{cc} \hat U_1 & \hat U_2\rix \mat{cc}\hat \Sigma_1 & 0\\ 0 & 0 \rix \mat{c}\hat U_1^T \\ \hat U_2^T \rix$ be a SVD of $\hat H_{11}$ with $\hat U_1 \in \R^{r,s}$ and $\hat \Sigma_1 \in \R^{s,s}$. Define 
  $\hat H_{11}^\epsilon=\mat{cc} \hat U_1 & \hat U_2\rix \mat{cc}\hat \Sigma_1 & 0\\ 0 & \Gamma \rix \mat{c}\hat U_1^T \\ \hat U_2^T \rix$, where $\Gamma$ is a diagonal matrix with diagonal entries each equal to $\frac{\epsilon}{\beta}$, where $\beta$ is given by~\eqref{eq:def_beta}, and set $\hat A_{11}^\epsilon =\hat H_{11}^\epsilon+\hat S_{11}$. Clearly $\hat A_{11}^\epsilon \in \mathcal{N}_{\succeq}^{S_r} $, and 
\begin{eqnarray}\label{eq:proof3}
{\|\hat{A}_{11}^{\epsilon} \Sigma_1- U_1^T B V_1\|}_F^2 &=& 
{\|(\hat{H}_{11}^{\epsilon}+\hat S_{11}) \Sigma_1- U_1^T B V_1\|}_F^2 \nonumber\\
&=& {\|(\hat{H}_{11}+\hat S_{11}) \Sigma_1- U_1^T B V_1+ \hat U_2 \Gamma \hat U_2^T \Sigma_1\|}_F^2 \nonumber \\
&\leq &  {\|\hat{A}_{11} \Sigma_1- U_1^T B V_1\|}_F^2 + {\|\Gamma\|}_F^2 {\|\Sigma_1\|}^2_F + 2 {\|\hat{A}_{11} \Sigma_1- U_1^T B V_1\|}_F {\|\Gamma\|}_F {\|\Sigma_1\|}_F \nonumber \\
&< &   {\|\hat{A}_{11} \Sigma_1- U_1^T B V_1\|}_F^2  + \epsilon,
\end{eqnarray} 
 where the last inequality follows since $\epsilon$ satisfies \eqref{epsilon_inequality_eq} and ${\|\Gamma\|}_F= \frac{\epsilon}{4{\|\Sigma_1\|}_F {\|\hat{A}_{11} \Sigma_1- U_1^T B V_1\|}_F}$ when ${\|\hat{A}_{11} \Sigma_1- U_1^T B V_1\|}_F \neq 0$ and $\frac{\epsilon}{4{\|\Sigma_1\|}_F}$ otherwise. 
Now define	 
\begin{equation}
A_{\epsilon} := U \left( \begin{bmatrix}\hat{H}_{11}^\epsilon & \frac{1}{2}(U_2^T BV_1 \Sigma_1^{-1})^T\\ \frac{1}{2}U_2^T BV_1 \Sigma_1^{-1} & K_\epsilon \end{bmatrix}
+ \begin{bmatrix} \hat{S}_{11} & -\frac{1}{2}(U_2^T BV_1 \Sigma_1^{-1})^T\\ \frac{1}{2}U_2^T BV_1 \Sigma_1^{-1} & R \end{bmatrix} \right)U^T
\end{equation} 
for some $K_{\epsilon} \in \mathbb{R}^{n-r,n-r}$ such that 
$K_{\epsilon}-\frac{1}{4}(U_2^{T} B V_1 \Sigma_1 ^{-1}) (\hat{H}_{11}^{\epsilon})^{-1} (U_2 ^T B V_1 \Sigma_1 ^{-1})^T \succeq 0$ and $R \in \R^{n-r,n-r}$ such that $R^T=-R$.  This yields \eqref{A_eps} after simplifications,
 and by construction $A_{\epsilon} \in \mathcal{N}_{\succeq}^{S_n}$.
Thus from~\eqref{eq:proof1} we have 
\begin{eqnarray} \label{A_eps_case2}
 {  \|A_{\epsilon}X-B\|}_F^2 &=& {\|(\hat{H}_{11}^{\epsilon}+\hat S_{11)} \Sigma_1- U_1^TBV_1\|}_F^2 + {\|BV_2\|}_F^2\nonumber \\
& =& {\|\hat{A}_{11}^\epsilon \Sigma_1- U_1^TBV_1\|}_F^2 + {\|BV_2\|}_F^2 \nonumber \\
&<& {\|\hat{A}_{11} \Sigma_1- U_1^TBV_1\|}_F^2 + {\|BV_2\|}_F^2 + \epsilon,
 \end{eqnarray}
where the last identity follows from~\eqref{eq:proof3}.
As $\epsilon$ tends to zero, from~\eqref{eq:proof2} and~\eqref{A_eps_case2}, we get the equality in~\eqref{eq:proof2}.
Hence 
\begin{equation}
 \inf_{A \in \mathcal{N}_{\succeq}^{S_n}} {  \|AX-B\|}_F^2 = {\|\hat{A}_{11} \Sigma_1- U_1^TBV_1\|}_F^2 + {\|BV_2\|}_F^2.
 \end{equation}
This infimum is attained when $\epsilon =0$, however from
Lemma~\ref{lem:psd_character} when $\epsilon =0$,  $A_{\epsilon} \not\in \mathcal{N}_{\succeq}^{S_n}$ as
$\text{null}(\hat{H}_{11}) \not \subseteq \text{null}(\frac{1}{2}U_2^T BV_1 \Sigma_1^{-1})$.  Therefore by the fact that $\inf_{K \in \R^{n-r,r}}{\| K\Sigma_1-U_2^TBV_1\|}_F^2=0$ and the uniqueness of $\hat A_{11}$ imply that the infimum is not attained. This completes the proof.
\eproof			  

Next, we obtain a corollary that characterizes the minimal Frobenius norm solutions to problem~\eqref{nspsdproc}.
 \begin{corollary}
Let $X,B \in \mathbb{R}^{n,m},$ and let $r=\rank(X).$ Let also $U_1,U_2,V_1,V_2,\Sigma_1,\hat{A}_{11}$, and $\hat H_{11}$  be as defined in Theorem~\ref{thetheorem}, and $Z:=\frac{1}{2}U_2^T B V_1 \Sigma_1^{-1}.$
\begin{enumerate}
\item If $\text{{\rm null}}(\hat{H}_{11}) \subseteq \text{{\rm null}}(Z)$, then $A_{opt}$ in~\eqref{A_opt} with $K=Z\hat{H}^{\dagger}_{11} Z^T$ and $R=0$,  is the  unique solution of the problem~\eqref{nspsdproc} with minimal Frobenius norm, that is, 
\begin{equation}
\left\{ A_{opt} \right\} = {\rm argmin}_{A \in{\rm argmin}_{A \in \mathcal{N}_{\succeq}^{S_n}} {\|AX-B\|}_F} {\|A\|}_{F}.
\end{equation} 
\item Otherwise, for sufficiently small $\epsilon > 0$ according to \eqref{epsilon_inequality_eq}, the matrix $A_{\epsilon}$ in \eqref{A_eps} with $K_{\epsilon}= Z (\hat{H}^{\epsilon}_{11})^{-1} Z^T$  and $R=0$, is the unique matrix in $\mathcal{N}_{\succeq}^{S_n}$ with minimal Frobenius norm
 such that 
\begin{equation}
{\|A_{\epsilon}X-B\|}_F^2 < \inf_{A \in \mathcal{N}_{\succeq}^{S_n}} {\|AX-B\|}_F^2 + \epsilon.
\end{equation}
\end{enumerate}
 \end{corollary}
 \proof In view of Lemma~\ref{lem:minnorm}, the proof is similar to~\cite[Corollary 1]{GilS18}. 
 \eproof

Note that when $\text{rank}(X)=m$ and $BX^T+XB^T \preceq 0$, then the PSDP problem has a unique solution $A=0$~\cite[Theorem 2.5]{Woo96}. In the following, we show that a similar result is also true for NSPSD matrices. In particular, when $U_1^T(BX^T+XB^T) \preceq 0$ (where $U_1$ is as defined in Theorem~\ref{thetheorem}) then the subproblem~\eqref{semi_analytic_eq} has a unique solution $\hat A_{11}=0$ and as a result the exact value of the infimum in~\eqref{NSPSD Problem} can be computed. This is stated in the following theorem proof of which is similar to~\cite[Theorem 2]{GilS18}. 
\begin{theorem}
Let $X,B \in \mathbb{R}^{n,m},$ and let $r=\rank(X)<n.$ Let also $U_1,U_2,V_1,V_2$ and $\Sigma_1$ be as defined in Theorem~\ref{thetheorem}. If $U_1^T(BX^T + XB^T)U_1 \preceq 0,$ then 
\begin{equation}
\inf_{A \in \mathcal{N}_{\succeq}^{S_n}} {\|AX-B\|}_F^2 = {\|U_1^T B V_1 \|}_F^2 + {\|BV_2\|}_F^2,
\end{equation}
and it is not attained for any $A \in \mathcal{N}_{\succeq}^{S_n}.$ In this case, let $\epsilon >0$ be sufficiently small and let $A^{\epsilon}_{11} \in \mathbb{R}^{r,r}$ be a diagonal matrix with diagonal entries each equal to $\frac{\epsilon}{\alpha},$ where $\alpha=4 \sqrt{n} {\|\Sigma_1\|}_F {\|U_1^T BV_1\|}_F.$ Define
\begin{equation}
A_{\epsilon}:= U_1 A_{11}^{\epsilon}U_1^T + U_2(U_2^T BV_1 \Sigma_1^{-1}) U_1^T +  U_2K^{\epsilon} U_2^{T}+ U_2RU_2^{T},
\end{equation}
where $K_{\epsilon} \in \mathbb{R}^{n-r,n-r}$ is such that $K_{\epsilon}-\frac{1}{4}(U_2^T BV_1 \Sigma_1^{-1})(H^{\epsilon}_{11})^{-1} (U_2^T BV_1 \Sigma_1^{-1})^{T} \succeq 0 $, and $R \in \R^{n-r,n-r}$ such that $R^T=-R$. Then $A_{\epsilon} \in \mathcal{N}_{\succeq}^{S_n}$ and 
\begin{equation}
{\|A_{\epsilon}X-B\|}_F^2 < \inf_{A \in \mathcal{N}^{S_n}_{\succeq} } {\|AX-B\|}_F^2 + \epsilon.
\end{equation}
\end{theorem}
\proof In view of Lemma~\ref{lem:trace}, the proof is similar to~\cite[Theorem 2]{GilS18}.
\eproof		  
A particular case of Theorem~\ref{thetheorem} is the vector case (i.e., when $m=1$) or more generally when $\text{rank}(X)=1$, where the assumption on the solution to the subproblem~\eqref{assumtion} is not required. In the following, we obtain a result for the case $\text{rank}(X)=1$ that solves the NSPSDP problem~\eqref{nspsdproc} exactly whenever the infimum is attained. This is analogus to~\cite[Theorem 3]{GilS18} where it was stated for the PSDP problem. 
\begin{theorem}
Let $X,B \in \mathbb{R}^{n,m}$ be such that $\rank(X)=1$. Let 
$X=U \Sigma V^{T}$ be a SVD of $X$, where $U=[u ~ U_1] \in \mathbb{R}^{n,n}$ 
	with $u \in \mathbb{R}^{n},V =[v ~ V_1] \in \mathbb{R}^{m,m}$ with $v \in \mathbb{R}^{m},$ and $\Sigma= \begin{bmatrix}\sigma & 0\\0 & 0\end{bmatrix} \in \mathbb{R}^{n,m}$ with $\sigma >0$. Then the following hold.
\begin{enumerate}
\item If $u^T B v >0$, then 
\begin{equation}
\inf_{A \in \mathcal{N}_{\succeq}^{S_n} } {\|AX-B\|}_F = {\|BV_1\|}_F,
\end{equation}
and $A_{opt}$ attains the infimum if and only if 
\begin{equation} \label{A_opt_vec}
A_{opt}= \sigma^{-1}((u^T B v)uu^T + U_1 U_1 ^T B vu^T ) + U_1K U_1 ^T+U_1RU_1^T,
\end{equation}
for some matrix $K\in \R^{n-1,n-1}$  such that $K- \frac{1}{4\sigma  u^T B v } U_1^T (Bv) (Bv)^T U_1$ and $R \in \R^{n-1,n-1}$ such that $R^T=-R$.
\item If $u^T B v \leq 0$, then 
\begin{equation} \label{semi_analytic_eq_vec}
\inf_{A \in \mathcal{N}_{\succeq}^{S_n}} {\|AX-B\|}_F^2 = {\|u^T B v\|}_F^2 + {\|BV_1\|}_F^2.
\end{equation}
Further, if $U_1 ^T B v =0$, then the infimum in \eqref{semi_analytic_eq_vec} is attained by a matrix $\widetilde A_{opt}$ of the form 
\[
\widetilde A_{opt}=U_1K U_1 ^T+U_1RU_1^T,
\]
where $K,R\in \R^{n-1,n-1}$  such that $K\succeq 0$ and $R^T=-R$.
If $U_1^T Bv \neq 0$, then the infimum in \eqref{semi_analytic_eq_vec} is not attained. In the later case for any arbitrary small $\epsilon > 0$, choose $n_0 \in \mathbb{N}$ such that $\frac{\sigma^2}{n_0 ^2}- 2 \sigma \frac{u^T Bv}{n_0} < \epsilon$ and define 
\begin{equation}
A_{n_0}= \sigma ^{-1} ( \frac{1}{n_0} u u^T +U_1U_1 ^T B v u^T ) + U_1 K_{n_0} U_1^T+U_1RU_1^T,
\end{equation}
for some $K_{n_0}$ with $K_{n_0}- \frac{n_0}{4 \sigma} U_1^T (B v)(B v)^T U_1 \succeq 0$. Then $A_{n_0} \in \mathcal{N}_{\succeq}^{S_n}$ and 
\begin{equation}
{\|A_{n_0}X-B\|}_F^2  < \inf_{A \in \mathcal{N}_{\succeq}^{S_n}}{\|AX-B\|}_F^2 + \epsilon.
\end{equation}
\end{enumerate} 
\end{theorem}

\section{On the complex NHPSD matrix Procrustes problem}\label{sec:compnspsdp}

In this section, we consider the complex non-Hermitian positive semidefinite procrustes (NHPSDP) problem~\eqref{complex_main_problem}. 
Motivated by the work in~\cite{KisH07} for complex Hermitian positive semidefinite matrices, we show that the complex NHPSDP problem can be equivalently transformed in a real NSPSDP problem. If $X$ is of full column rank, then this formulation results in a real NSPSDP problem of size double the original problem, for which the method developed in Section~\ref{sec:algo} is applicable. 
On the other hand, if $X$ is column-rank deficient, then this formulation results  in a real structured NSPSDP problem, where the NSPSD matrices belong to the set $ \mathcal{N_S}_\succeq ^{S_{2n}}$ defined by 
\begin{equation}\label{def:struns}
 \mathcal{N_S}_\succeq ^{S_{2n}}:=\left\{ P=\mat{cc} P_1 & P_2 \\ -P_2 & P_1 \rix \in \mathcal{N}_\succeq ^{S_{2n}}:~P_1,\,P_2 \in \R^{n,n}
 \right\}.
\end{equation}
This formulation is stated in the following result without its proof  which is similar to~\cite[Theorem 3.3]{KisH07}. 

\begin{theorem}\label{thm:main_comp}
Let {{$X,B \in \mathbb{C}^{n,m} \setminus \{0\}$}}, where $n\geq m$,  and let $X=X_r + i X_j$  and $B=B_r + i B_j$, where $X_r,X_j,B_r,B_j \in \mathbb{R}^{n,m}$. Consider the three optimization problems,
\begin{equation}\label{complex_nspsd}
\inf_{A \in \mathcal{N}_\succeq^{H_n}} {\| AX-B\|}_F, \tag{$\mathcal{P}_\C$}
\end{equation} 
 \begin{equation}\label{real_strnspsd}
\inf_{P \in \mathcal{N_S}_\succeq ^{S_{2n}}} \begin{Vmatrix}
P \begin{pmatrix}X_r & X_j\\-X_j & X_r\end{pmatrix} -
\begin{pmatrix}B_r & Bj\\-B_j & B_r\end{pmatrix}\end{Vmatrix},
\tag{$\mathcal{P'}$}
\end{equation}
and 
\begin{equation}\label{real_nspsd}
\inf_{\widetilde A \in \mathcal{N}_\succeq ^{S_{2n}}} \begin{Vmatrix}
\widetilde A \begin{pmatrix}X_r & X_j\\-X_j & X_r\end{pmatrix} -
\begin{pmatrix}B_r & Bj\\-B_j & B_r\end{pmatrix}\end{Vmatrix}.\tag{$\mathcal{P''}$}
\end{equation}
Then 
\begin{enumerate}
\item If $\text{\rm rank}(X)=m$, then~\eqref{real_nspsd}=~\eqref{real_strnspsd}=$\sqrt{2}$\,\eqref{complex_nspsd}. Moreover, the infimum in~\eqref{complex_nspsd} and~\eqref{real_nspsd} is attained for some unique 
$A_0 \in \mathcal{N}_\succeq^{H_n} $ and $\widetilde A_0 \in \mathcal{N}_\succeq ^{S_{2n}}$, respectively. In this case, write
 $\widetilde A_0 = \begin{bmatrix} A_1 & A_2\\ A_3 & A_4\end{bmatrix}$,
where $A_1,A_2,A_3,A_4 \in \mathbb{R}^{n,n}$. Then $A_1=A_4$, $A_3=-A_2$, and $A_0=A_1 + i A_2$. 
\item  If $\text{\rm rank}(X)<m$, then~\eqref{real_nspsd} $<$~\eqref{real_strnspsd} $=\sqrt{2}$\,\eqref{complex_nspsd}. 
\end{enumerate}	
\end{theorem}
\proof In view of Lemma~\ref{complex_real_nspsd_lemma} and Theorem~\ref{lem:KrisNSPSDP} , the proof can be derived similar to~\cite[Theorem 3.3]{KisH07}.
\eproof
Note that, when $X$ is not of full column rank then from Theorem~\ref{thm:main_comp}, a solution to~\eqref{real_nspsd} will give a lower bound to the infimum in~\eqref{complex_nspsd}. On the other hand, if we take into consideration the structure of the set $\mathcal{N_S}_\succeq ^{S_{2n}}$ in~\eqref{real_strnspsd}, then a projected FGM can be used for solving~\eqref{real_strnspsd}. For this we need to calculate the projection of a block $2n \times 2n$ matrix onto the set 
$\mathcal{N_S}_\succeq ^{S_{2n}}$ of block non-symmetric semidefinite matrices defined in~\eqref{def:struns}. This is achieved in the following result, a proof of which is provided in~\ref{subsec:proof}. 
Recall that any complex matrix $A \in \C^{n,n}$ can be denoted by  $A=A_1+iA_2$, where $A_1,A_2 \in \R^{n,n}$, and $R(A)=\mat{cc}A_1 &A_2 \\-A_2 & A_1\rix$.

\begin{theorem}\label{thm:complex_nhpsdproj}
Let $A=\mat{cc} A_1&A_2 \\ A_3 & A_4\rix$, where $A_1,A_2,A_3,A_4 \in \R^{n,n}$. Then 
\begin{equation}
{\rm argmin}_{X \in \C^{n,n},\, R(X) \in \mathcal{N_S}_\succeq ^{S_{n}} } {\left \|A-R(X) \right \|}_F =
\mat{cc}  {\widetilde X_1}{}_H+ {\widetilde {X}_1}{}_S &
{\widetilde X_2}{}_S+{\widetilde X_2}{}_H \\ -{\widetilde X_2}{}_S-{\widetilde X_2}{}_H  & {\widetilde X_1}{}_H+{\widetilde X_1}{}_S \rix,
\end{equation}
where ${\widetilde X_1}{}_H+ i{\widetilde {X}_2}{}_S $ is the PSD projection of 
$\frac{1}{2}({A_1}_H+{A_4}_H)+i\frac{1}{2}({A_2}_S-{A_3}_S)$, 
${\widetilde X_1}{}_S=\frac{{A_1}_S+{A_4}_S}{2}$, and 
${\widetilde X_2}{}_H=\frac{{A_2}_H-{A_3}_H}{2}$.
\end{theorem}

\section{Algorithm for the NSPSDP problem} \label{sec:algo}

Building on the theoretical results presented in the previous sections, we now propose a method 
to solve the NSPSDP problem~\eqref{nspsdproc}. 
Our proposed algorithm uses the same strategy as in~\cite{GilS18} for the symmetric PSDP problem, which was shown to be very efficient. The two main differences with the algorithm proposed in~\cite{GilS18} are 
(1) the semi-analytical approach to reduce the problem size (see Section~\ref{sec:nspsdp}), 
and  
(2) the projection step in our iterative algorithm, which is  performed onto the set of NSPSD matrices instead of the PSD matrices.  
Hence we do not dig deep into the technical derivations but rather provide the high level ideas. We refer the interested reader to~\cite{GilS18}  for more details.

This algorithm consists in two main steps: 
\begin{enumerate}

\item It uses the semi-analytical solution described in Section~\ref{sec:nspsdp} to reduce the NSPSDP problem with $X,B \in \mathbb{R}^{n \times m}$ to a diagonal NSPSDP problem with $\tilde X, \tilde B \in \mathbb{R}^{r \times r}$ where $r = \rank(X)$ and $\tilde X$ is diagonal. 

\item It applies the fast gradient method (FGM) from~\cite{nes83} to the reduced diagonal NHPSDP problem. 
FGM is an optimal first-order method, and is guaranteed to decrease the objective function values at a linear rate $O( (1-1/\kappa)^t )$  for strongly convex problems, where $t$ is the iteration count and $\kappa = \frac{\sigma_1(X)}{\sigma_r(X)} > 0$. This is much faster than the standard gradient descent method, with a linear rate of $O( (1-1/\kappa^{2})^t )$; see\cite{nes04} for more details.

\end{enumerate}

\begin{remark} {\rm The semi-analytical solution that reduces the problem size and makes $X$ diagonal can be combined with any other method. We use here the FGM because it scales well, being a first-order method. Moreover, second-order methods, in particular interior-point methods, have been proposed for this problem; see Section~\ref{sec:numerical}.  
}
\end{remark}

The last ingredient of the proposed method is the choice of the initial matrix $A$ to start the iterative algorithm. Again, we follow the strategy proposed in~\cite{GilS18}. It is based on the following two observations: 
\begin{enumerate}

\item The convergence of FGM depends on the condition number of $X$ that is, the ratio of the largest to the smallest diagonal entries of $X$ since it is diagonal.   

\item The problem 
\[
\min_{A_{1} \succeq 0, A_{2} \succeq 0}
{\left\|
\left( \begin{array}{cc} A_{1} & 0 \\ 0 & A_{2} \end{array} \right)
\left( \begin{array}{cc} X_{1} & 0 \\ 0 & X_{2} \end{array} \right)
- \left( \begin{array}{cc} B_{1} & B_{12} \\ B_{21} & B_{2} \end{array} \right)   \right\|}_F , 
\]
 can be decoupled into two independent subproblems: for $i=1,2$,
\[
\min_{A_{i} \succeq 0} {\|A_{i} X_{i} - B_{i}\|}_F, 
\]
and this observation can be easily generalized to any number of blocks. 

\end{enumerate} 

The initialization in~\cite{GilS18} first partitions the diagonal of $X$ by entries of similar magnitude (namely, with ratio at most 100, using a recursive procedure) so that FGM applied on each subproblem converges linearly with rate at least 0.99, and 100 iterations of FGM  are applied on these subproblems to initialize $A$ as a block diagonal matrix.

%
%
%
%
%
%
%
%
%
%

\paragraph{Computational cost} The main computational cost of the semi-analytical approach is the SVD computation of $X$, which requires $\mathcal{O}(mn \min(m,n))$ operations. 

For the NSPSDP problem in dimension $n$, the FGM requires the computation of the gradient, $AXX^T-BX^T$, 
for a total of $\mathcal{O}(n^3 + mn^2)$ operations (note that $XX^T$ and $BX^T$ need to be computed only once).  Then it requires the projection onto the set of NSPSD matrices which requires an eigenvalue decomposition in  $\mathcal{O}(n^3)$ operations. 
If $X$ is diagonal with $m=n$, 
these costs reduce to $\mathcal{O}(n^3)$ operations. 
Using the semi-analytical approach reduces the cost of FGM from $\mathcal{O}(n^3 + mn^2)$ to $\mathcal{O}(r^3)$ operations. As we will see in the numerical experiments, this allows to significantly speed up FGM when $X$ is rank deficient. 

Finally, the computational cost of the proposed approach is 
$\mathcal{O}( mn \min(m,n) + T r^3 )$ where $T$ is the number of iterations of the FGM.

\section{Numerical experiments} \label{sec:numerical}

In this section, we compare our proposed algorithm combining a semi-analytic approach and the FGM, which we refer to as AN-FGM, to the following: 
\begin{itemize}


\item The IPM SDPT3 (version 4.0)~\cite{toh1999sdpt3, tutuncu2003solving}, where we used CVX as a modeling system~\cite{cvx, gb08}. We refer to this method as SDPT3. 

\item The FGM method 
without the semi-analytic approach, using the identity matrix as an initialization which is (optimally) scaled by the scalar  
\[
\argmin_{\alpha} \| \alpha X - B\|_F = \frac{\text{trace}(XB^T)}{ \|X\|_F^2}. 
\]

\end{itemize}

We first compare them on synthetic data sets, generated in the same way as in~\cite{GilS18}, and then to two data sets from~\cite{KrisLVP04} and~\cite{KisH07}. 

The codes and data sets are available from \url{http://bit.ly/NSPSDProc_v1}. Note that we have merged our proposed algorithm for the NSPSDP problem with the code from~\cite{GilS18}. This leads to a single Matlab file that can handle the symmetric and non-symmetric PSDP problems, in the real and complex cases. 

Interior point methods do not scale well as they need to solve an $n^2$-by-$n^2$ linear system at each iteration, with a cost of $\mathcal{O}(n^6)$ operations per iteration. 
In Matlab, they can handle matrices of size at most 100 within a few minutes, while FGM can handle much larger matrices, up to size 1000.

\begin{remark}{\rm 
We also tested with the two methods proposed in~\cite{KrisLVP04} that use  dedicated interior-point methods (IPMs); the code is available from \url{https://sites.google.com/site/nathankrislock/software}. 
However, these two methods were outperformed by SDPT3 in most cases, and hence we do not report their results here. 
}
\end{remark}

\subsection{Synthetic data sets from~\cite{GilS18}} 

We first use the synthetic data sets from~\cite{GilS18}.  We have $m=n$, $m < n = 2m$ or $n < m=2n$ and, in all cases, we set $\max(m,n) = 60$. 
The matrix $B$ is generated in the same way: each entry is randomly generated following a normal distribution with mean 0 and standard deviation 1 (\texttt{randn(m,n)} in Matlab).
For the matrix $X$, there are three cases
\begin{enumerate}
\item Well-conditioned. Each entry is randomly generated following a normal distribution with mean 0 and standard deviation 1 (\texttt{randn(m,n)} in Matlab).
\item Ill-conditioned. Let $(U,\Sigma,V)$ be the compact SVD of a matrix generated as in the well-conditioned case. Then we generate $X = U \Lambda V$, where $\Lambda$ is a diagonal matrix such that $\Lambda(i,i) = \alpha^{i-1}$ and $\alpha^{\min(m,n)-1} = 10^6 = \kappa(X)$.
\item Rank deficient. We perform the SVD $(U,\Sigma,V)$ of a matrix generated as in the well-conditioned case, set the $r = \min(m,n)/2$ smallest singular values
of $\Sigma$ to zero to obtain $\Sigma'$, and then compute $X = U \Sigma' V^T$ so that $\rank(X) = \min(m,n)/2$.
\end{enumerate}
There are therefore a total of 9 scenarios (3 cases for the choice of $n$ and $m$, and three choices for the generation of $X$). 
We run all algorithms with their default stopping criterion. 
We set the maximum run time of FGM as the time taken by IPM-CVX, while we also stop the algorithm if an iterate $A^{(k)}$ satisfies 
\[
\| A^{(k)} - A^{(k-1)} \|_F < \delta \| A^{(1)}  -  A^{(0)} \|_F
\]
 using $\delta = 10^{-6}$, that is, the modification of $A$ compared to the first step is less than  $10^{-6}$. 

Table~\ref{tab:results} reports the average and standard deviation of the final relative error (in percent), defined as $\|AX-B\|_F / \|B\|_F$, as well as the run time for each algorithm, over 20 randomly generated matrices for type of synthetic data set. 
 \begin{center}  
 \begin{table}[h!] 
 \begin{center} 
\caption{Average relative error (rel.\ err., in \%) and runtime (s.), as well as their standard deviations
for CVX, FGM, AN-FGM on 20 randomly generated problems for 9 different scenarios. 
The best result is highlighted in bold. \label{tab:results} } 
 \begin{tabular}{|c|c|c||c|c|c|} 
 \hline          &        &     &        SDPT3        &        FGM       &       AN-FGM      \\ \hline  \hline  
 Well-cond. & $m$=$n$=60 & rel. err. (\%)  &   \textbf{18.37 $\pm$ 0.67}    &  \textbf{18.37 $\pm$ 0.67}  &   \textbf{18.37  $\pm$ 0.67}  \\ 
                 &          & time (s.)   &   21.37 $\pm$ 2.41    &  5.49 $\pm$ 4.98   &   \textbf{1.43  $\pm$ 0.26}  \\ \hline 
                 & $m$=$2n$=$60$ & rel. err. (\%)   &   \textbf{26.56 $\pm$ 0.97}    &  \textbf{26.56 $\pm$ 0.97}  &   \textbf{26.56  $\pm$ 0.97}  \\ 
                 &          & time (s.)   &   3.09 $\pm$ 0.56     &  \textbf{0.03 $\pm$ 0.02}   &   0.09  $\pm$ 0.03  \\ \hline 
                 & $n$=$2m$=$60$ & rel. err. (\%)   &   \textbf{17.43 $\pm$ 0.94}    &  \textbf{17.43 $\pm$ 0.94}  &   \textbf{17.43  $\pm$ 0.94}  \\ 
                 &          & time (s.)   &   21.09 $\pm$ 2.77    &  5.19 $\pm$ 2.54   &   \textbf{0.10  $\pm$ 0.03}  \\ \hline \hline 
 Ill-cond. & $m$=$n$=$60$ & rel. err. (\%)    &   \textbf{19.41 $\pm$ 0.56}    &  20.31 $\pm$ 0.61  &   20.49  $\pm$ 0.73  \\ 
                 &          & time (s.)   &   34.69 $\pm$ 2.26    &  34.72 $\pm$ 2.26  &   \textbf{1.43  $\pm$ 0.13}  \\ \hline 
                 & $m$=$2n$=$60$ & rel. err. (\%)   &   \textbf{27.00 $\pm$ 0.85}    &  27.63 $\pm$ 0.62  &   27.71  $\pm$ 0.74  \\ 
                 &          & time (s.)   &   3.72 $\pm$ 0.43     &  3.75 $\pm$ 0.43   &   \textbf{0.44  $\pm$ 0.06}  \\ \hline 
                 & $n$=$2m$=$60$ & rel. err. (\%)   &   20.27 $\pm$ 2.87    &  24.02 $\pm$ 1.04  &   \textbf{20.16  $\pm$ 1.13}  \\ 
                 &          & time (s.)   &   50.58 $\pm$ 13.54    &  50.61 $\pm$ 13.55  &   \textbf{0.38  $\pm$ 0.07}  \\ \hline \hline 
 Rank-def.  & $m$=$n$=$60$ & rel. err. (\%)    &   \textbf{21.79 $\pm$ 0.74}    &  \textbf{21.79 $\pm$ 0.74}  &   \textbf{21.79  $\pm$ 0.74}  \\ 
                 &          & time (s.)   &   23.51 $\pm$ 13.24    &  3.34 $\pm$ 1.51   &   \textbf{0.05  $\pm$ 0.02}  \\ \hline 
                 & $m$=$2n$=$60$ & rel. err. (\%)   &   \textbf{27.57 $\pm$ 0.84}    &  \textbf{27.57 $\pm$ 0.84}  &   \textbf{27.57  $\pm$ 0.84}  \\ 
                 &          & time (s.)   &   2.95 $\pm$ 1.41     &  0.68 $\pm$ 0.30   &   \textbf{0.01  $\pm$ 0.02}  \\ \hline 
                 & $n$=$2m$=$60$ & rel. err. (\%)   &   \textbf{26.17 $\pm$ 0.64}    &  \textbf{26.17 $\pm$ 0.64}  &   \textbf{26.17  $\pm$ 0.64}  \\ 
                 &          & time (s.)   &   80.92 $\pm$ 111.00    &  2.87 $\pm$ 1.20   &   \textbf{0.00  $\pm$ 0.01}  \\ \hline 
\end{tabular} 
 \end{center} 
 \end{table} 
 \end{center} 
 
 We observe the following: 
 \begin{itemize}

\item For well-conditioned and rank-deficient problems, all algorithms find the same solution, although FGM-based algorithms run significantly faster than CVX, a second-order method based on IPM, which is expected. 

\item Except for well-conditioned case with $m=2n$ where the semi-analytical solution does not reduce the problem size (since $r=n$), AN-FGM runs significantly faster than FGM, especially for the rank-deficient cases. This follows from the fact that AN-FGM has a significantly lower computational cost than FGM when $r = \rank(X) \ll n$, and because the semi-analytical approach reduces the NSPSDP problem to a diagonal problem (hence strongly convex) with guaranteed linear convergence rate. 

\item Except for ill-conditioned case with $m=n=60$ and $m=2n=60$, AN-FGM outperforms CVX. The reason CVX performs better in these two cases is because  the linear rate of convergence of AN-FGM is rather slow, because of the ill-conditioned nature of the problem. Interestingly, for $n=2m=60$, AN-FGM provides more accurate solution, because the semi-analytical reduces the problem size, since $\rank(X) = m < n=2m$. 
 
 \end{itemize}
 
Note that similar observations were made in~\cite{GilS18} for the symmetric PSDP problem.

\subsection{Local compliance estimation problem from \cite{KrisLVP04}} 

In~\cite{KrisLVP04}, an NSPSDP problem is considered corresponding to a local compliance estimation problem, with the following matrices (with two digits of accuracy):  
\[
X^T = 
\left( \begin{array}{ccc} 
 -0.32 &  0.03 &  0.06 \\ 
 -0.33 &  -0.02 &  0.06 \\ 
 -0.36 &  0.08 &  0.06 \\ 
 -0.30 &  0.03 &  0.05 \\ 
 -0.32 &  -0.00 &  0.07 \\ 
 -0.34 &  0.07 &  0.05 \\ 
 -0.24 &  0.07 &  0.05 \\ 
 -0.21 &  -0.01 &  0.02 \\ 
 -0.33 &  0.16 &  0.10 \\ 
 -0.25 &  0.09 &  0.06 \\ 
 -0.22 &  0.00 &  0.03 \\ 
 -0.31 &  0.15 &  0.09 \\ 
\end{array} \right), 
B^T = 
\left( \begin{array}{ccc} 
 -1.43 &  0.15 &  -0.44 \\ 
 -1.40 &  -0.31 &  -0.42 \\ 
 -1.38 &  0.44 &  -0.42 \\ 
 -1.43 &  0.14 &  -0.44 \\ 
 -1.40 &  -0.31 &  -0.42 \\ 
 -1.37 &  0.43 &  -0.42 \\ 
 -1.43 &  0.16 &  -0.43 \\ 
 -1.40 &  -0.32 &  -0.42 \\ 
 -1.38 &  0.42 &  -0.43 \\ 
 -1.43 &  0.15 &  -0.44 \\ 
 -1.40 &  -0.33 &  -0.42 \\ 
 -1.37 &  0.42 &  -0.44 \\ 
\end{array} \right)  . 
\]
This is a small well-conditioned problem ($\kappa(X) = 28.1$), hence easy to solve. 
Using the same settings as in the previous section, all algorithm converge to the same solution, namely (with 4 digits of accuracy) 
\[
A = \left( \begin{array}{ccc} 
 5.0392 &  0.4423 &  1.5978 \\ 
 -0.6207 &  6.0223 &  -6.8559 \\ 
 1.8979 &  -0.4079 &  2.7600 \\ 
\end{array} \right)  . 
\] 
However, AN-FGM and FGM take about 0.05 seconds to compute this solution while CVX takes about 2 seconds.  
This is the same solution as obtained with the algorithms 
proposed in~\cite{KrisLVP04}, with relative error 
$\| AX - B \|_F / \| B \|_F  = 18.99\%$, and the eigenvalues of 
$A+A^T$ being 0, 
   10.28, and 17.36.

\subsection{Numerical example for the complex NSPSD problem} 

Using Theorems~\ref{thm:main_comp} and~\ref{thm:complex_nhpsdproj}, AN-FGM can readily be applied on the complex NSPSDP problem: first transform the complex problem to the real problem~\eqref{real_strnspsd} 
using Theorem~\ref{thm:main_comp}, and then replace the projection step within the FGM using Theorem~\ref{thm:complex_nhpsdproj}. 
Let us apply this algorithm to the example in~\cite{KisH07} with $X$ given by 
\[
\left( \begin{array}{cccc} 
 0.4694 &  0.5354 &  0.1326 &  -0.0787 \\ 
 -0.9036 &  0.5529 &  1.5929 &  -0.6817 \\ 
 0.0359 &  -0.2037 &  1.0184 &  -1.0246 \\ 
 -0.6275 &  -2.0543 &  -1.5804 &  -1.2344 \\ 
\end{array} \right) 
+ i \left( \begin{array}{cccc} 
 0.2888 &  -0.4650 &  -1.3573 &  -1.3813 \\ 
 -0.4293 &  0.3710 &  -1.0226 &  0.3155 \\ 
 0.0558 &  0.7283 &  1.0378 &  1.5532 \\ 
 -0.3679 &  2.1122 &  -0.3898 &  0.7079 \\ 
\end{array} \right) \] 
and $B$ by 
\[
\left( \begin{array}{cccc} 
 0.0112 &  -0.9898 &  1.1380 &  -0.3306 \\ 
 -0.6451 &  1.3396 &  -0.6841 &  -0.8436 \\ 
 0.8057 &  0.2895 &  -1.2919 &  0.4978 \\ 
 0.2316 &  1.4789 &  -0.0729 &  1.4885 \\ 
\end{array} \right) 
+ i 
\left( \begin{array}{cccc} 
 -0.5465 &  -0.8542 &  0.4853 &  -0.0793 \\ 
 -0.8468 &  -1.2013 &  -0.5955 &  1.5352 \\ 
 -0.2463 &  -0.1199 &  -0.1497 &  -0.6065 \\ 
 0.6630 &  -0.0653 &  -0.4348 &  -1.3474 \\ 
\end{array} \right). 
\] 
As reported in~\cite{KisH07}, the error for the symmetric PSDP problem is $\| AX - B \|_F = 4.19$, where $A \succeq 0$. 
Solving the corresponding NSPSDP problem, we obtain the same solution with the three algorithms, namely (with 2 digits of accuracy) 
\[
A = 
\left( \begin{array}{cccc} 
 0.55 &  -0.12 &  0.34 &  -0.16 \\ 
 -0.59 &  0.82 &  -0.02 &  -0.33 \\ 
 0.03 &  -0.42 &  0.08 &  0.22 \\ 
 0.42 &  -0.17 &  -0.07 &  0.08 \\ 
\end{array} \right)  
+ i 
\left( \begin{array}{cccc} 
 -0.07 &  0.24 &  -0.14 &  0.31 \\ 
 -0.89 &  -0.57 &  -0.17 &  -0.16 \\ 
 0.05 &  0.09 &  -0.17 &  0.04 \\ 
 0.62 &  -0.10 &  0.10 &  -0.16 \\ 
\end{array} \right), 
\]
with $\| AX - B \|_F = 3.04$ which is, as expected, smaller than the symmetric PSDP solution. The solution $A$ is such that $A+A^T$ has a single non-zero eigenvalue equal to 3.04. 
Again, this is a small well-conditioned problem (in the transformed real problem, $\kappa(X) = 5.3$), hence easy to solve. 
However, AN-FGM and FGM take 0.03 seconds to compute this solution while CVX takes about 2 seconds.

\section{Conclusion} 

In this paper, we proposed a semi-analytical algorithm for the NSPSDP problem which solves a reduced and well-posed subproblem using a fast gradient method.  
More precisely, the problem is first reduced to a smaller NSPSDP problem that always has a unique solution: assuming the solution to the smaller problem is known, we can obtain an approximate solution of the original problem analytically for any accuracy; see Theorem~\ref{thetheorem}. 
To solve the subproblem, we have developed an efficient first-order method. 
Since the complex problem can be equivalently rewritten as an overparametrized real problem (Theorem~\ref{thm:main_comp}), 
our algorithm can also be used in this setting. 
We illustrated the efficiency of the proposed algorithm on several numerical examples. 


\section*{Acknowledgement} 

The author are thankful to Nathan Krislock for providing them with the numerical example from his paper~\cite{KrisLVP04}.

\bibliographystyle{siam}
\bibliography{bibliostable}

\newpage 

\appendix
\section{ Computing projections }


Recall that for a complex matrix $X=X_1+iX_2 \in \C^{n,n}$, where $X_1,X_2 \in \R^{n,n}$, $R(X)$ is defined by $R(X)=\mat{cc}X_1 &X_2 \\-X_2 & X_1\rix$. 
Let $A$ be a block matrix of the form $A=\mat{cc} A_1&A_2 \\ A_3 & A_4\rix$, where $A_1,A_2,A_3,A_4 \in \R^{n,n}$. In this section, we compute three different projections of $A$. The first one is the projection of $A$ onto the set of block matrices of the form $R(X)$, the second one is the projection of $A$ onto the symmetric positive semidefinite block matrices of the form $R(X)$, and the third one is the projection of $A$ onto the set of non-symmetric positive semidefinite block matrices of the form $R(X)$. The first two projections are obtained in~\ref{subsec:projresher} and the third one (proof of Theorem~\ref{thm:complex_nhpsdproj}) is obtained in~\ref{subsec:proof}.

\subsection{ Projection for Hermitian PSD matrices}\label{subsec:projresher}

%
%

The following result computes the nearest matrix of the form $R(X)$ to a given block matrix $A$. 

\begin{theorem}
Let $A=\mat{cc} A_1&A_2 \\ A_3 & A_4\rix$, where $A_1,A_2,A_3,A_4 \in \R^{n,n}$. Then 
\begin{equation}
{\rm argmin}_{X \in \C^{n,n}} {\left \|A-R(X) \right \|}_F =R(\widetilde {X}_1+i\widetilde {X}_2),
\end{equation}
where $\widetilde X_1=\frac{A_1+A_4}{2}$ and $\widetilde X_2=\frac{A_2-A_3}{2}$.
\end{theorem}
\proof
We have,
\begin{eqnarray}\label{eq:prothm11}
{\rm argmin}_{X \in \C^{n,n}} {\left \| A-R(X) \right \|}_F^2 &=&
{\rm argmin}_{X_1,X_2 \in \R^{n,n}}{\left\|\mat{cc} A_1&A_2 \\ A_3 & A_4\rix-\mat{cc} X_1&X_2 \\ -X_2 & X_1\rix\right \|}_F^2 \nonumber \\
&=& {\rm argmin}_{X_1 \in \R^{n,n}} \left( {\| A_1-X_1  \|}_F^2 +{\| A_4-X_1  \|}_F^2 \right) + \nonumber\\ 
 &\hspace{1cm} & {\rm argmin}_{X_2 \in \R^{n,n}} \left( {\| A_2-X_2  \|}_F^2 +{\| A_3+X_2  \|}_F^2 \right).
\end{eqnarray}
The first infimum in~\eqref{eq:prothm11}
is attained by $\widetilde X_1 =\frac{A_1+A_4}{2}$ and the second infimum is
attained by $\widetilde X_2 =\frac{A_2-A_3}{2}$. Using this in~\eqref{eq:prothm11} yields that 
\[
{\rm argmin}_{X \in \C^{n,n}} {\left \|A-R(X) \right \|}_F =R(\widetilde {X}), \quad \text{where}~\, \widetilde X=\widetilde {X}_1+i\widetilde {X}_2.
\]
\eproof

The following result computes the nearest symmetric PSD matrix of the form $R(X)$ to a given block matrix $A$. This projection can be used to find an approximate solution to the complex PSD Procrustes problem~\cite{KisH07}.

\begin{theorem}\label{thm:strucsymproj}
Let $A=\mat{cc} A_1&A_2 \\ A_3 & A_4\rix$, where $A_1,A_2,A_3,A_4 \in \R^{n,n}$. Then 
\begin{equation}
{\rm argmin}_{X \in \C^{n,n},\, X \succeq 0} {\left \|A-R(X) \right \|}_F =R(\widetilde X),
\end{equation}
where $\widetilde X=\widetilde X_1+i\widetilde X_2$ is the PSD projection of $A$, i.e., $\widetilde X=\mathcal P_{\succeq }(\widetilde A_1 +i\widetilde A_2)$ with $\widetilde A_1= \frac{{A_1}_H+ {A_4}_H}{2}$ and $\widetilde A_2= \frac{{A_2}_S- {A_3}_S}{2}$.
\end{theorem}
\proof 
For any $X=X_1+iX_2$, where  $X_1,X_2 \in \R^{n,n}$ such that $X$ is Hermitian PSD, we have $X_1^T=X_1$ and $X_2^T=-X_2$. This implies that 
\begin{eqnarray}\label{eq:projthm21}
{\|A-R(X)\|}_F^2 &=& {\left\| \mat{cc} A_1&A_2 \\ A_3 & A_4\rix -
\mat{cc} X_1&X_2 \\ -X_2 & X_1\rix \right \|}_F^2 \nonumber \\
&=& {\left\| \mat{cc}  {A_1}_H & \frac{A_2+A_3^T}{2} \\ \frac{A_3+A_2^T}{2} & {A_4}_H \rix + \mat{cc}  {A_1}_S & \frac{A_2-A_3^T}{2} \\ \frac{A_3-A_2^T}{2} & {A_4}_S \rix
-\mat{cc} X_1&X_2 \\ -X_2 & X_1\rix \right \|}_F^2 \nonumber \\
&=& {\left\| \mat{cc}  {A_1}_H & \frac{A_2+A_3^T}{2} \\ \frac{A_3+A_2^T}{2} & {A_4}_H \rix 
-\mat{cc} X_1&X_2 \\ -X_2 & X_1\rix \right \|}_F^2 + {\left \| \mat{cc}  {A_1}_S & \frac{A_2-A_3^T}{2} \\ \frac{A_3-A_2^T}{2} & {A_4}_S \rix \right \|}_F^2, 
\end{eqnarray}
where the last identity holds since for any matrix $Z \in \R^{n,n}$, we have 
$\|Z\|_F^2=\|Z_H\|_F^2+\|Z_S\|_F^2$. Thus by taking the infimum 
in~\eqref{eq:projthm21} over all complex PSD matrices, we have 
\begin{equation}\label{eq:projthm22}
{\rm argmin}_{X \in \C^{n,n},X \succeq 0} {\left \|A-R(X) \right \|}_F =
{\rm argmin}_{X=X_1+iX_2 \in \C^{n,n},\, X \succeq 0}{\left\| \mat{cc}  {A_1}_H & \frac{A_2+A_3^T}{2} \\ \frac{A_3+A_2^T}{2} & {A_4}_H \rix 
-\mat{cc} X_1&X_2 \\ -X_2 & X_1\rix \right \|}_F^2.
\end{equation}
By setting $U=\frac{1}{\sqrt{2}}\mat{cc}I_n &iI_n \\iI_n & I_n\rix$ unitary and using the fact that $\|\cdot\|_F$ is unitarily invariant, we have 
from~\eqref{eq:projthm22} that 
\begin{eqnarray}\label{eq:projthm23}
&& \hspace{-1.7cm} {\rm argmin}_{X \in \C^{n,n},X \succeq 0} {\left \|A-R(X) \right \|}_F  \nonumber \\
&=& {\rm argmin}_{X=X_1+iX_2 \in \C^{n,n},\, X \succeq 0}{\left\| U^*\mat{cc}  {A_1}_H & \frac{A_2+A_3^T}{2} \\ \frac{A_3+A_2^T}{2} & {A_4}_H \rix U
-U^*\mat{cc} X_1&X_2 \\ -X_2 & X_1\rix U \right \|}_F^2 \nonumber \\
&=& {\rm argmin}_{X=X_1+iX_2 \in \C^{n,n},\, X \succeq 0} 
{\left \| \mat{cc} \widetilde A_1 + i \widetilde A_2 & K \\ \overline K & \widetilde A_1 - i \widetilde A_2 
\rix -\mat{cc}X &0 \\ 0 & \overline X \rix\right \|}_F^2, 
\end{eqnarray}
where $\widetilde A_1= \frac{{A_1}_H+ {A_4}_H}{2}$, $\widetilde A_2= \frac{{A_2}_S- {A_3}_S}{2}$, and $K=(\frac{{A_3}_H+ {A_2}_H}{2})+ i ( \frac{{A_1}_H- {A_4}_H}{2})$. 
The last identity in~\eqref{eq:projthm23} is due to Lemma~\ref{lem:com1}.
 Thus, we have 
 \begin{eqnarray}\label{eq:projthm24}
 {\rm argmin}_{X \in \C^{n,n},X \succeq 0} {\left \|A-R(X) \right \|}_F^2
 &=&  {\rm argmin}_{X \in \C^{n,n},X \succeq 0} 
 {\left\| \widetilde A_1 + i \widetilde A_2  - X \right \|}_F^2 +
 {\left\| \overline{\widetilde A_1 + i \widetilde A_2}  - \overline X \right \|}_F^2 \nonumber \\
 &=& 2 \, {\rm argmin}_{X \in \C^{n,n},X \succeq 0} 
 {\left\| \widetilde A_1 + i \widetilde A_2  - X \right \|}_F^2.
 \end{eqnarray}
In view of Lemma~\ref{lem:com5}, the infimum in the right hand side of~\eqref{eq:projthm24} is attained by 
$P_{\succeq}(\widetilde A_1 + i \widetilde A_2)$, the PSD projection of $\widetilde A_1 + i \widetilde A_2 $ onto the set of complex PSD matrices. This concludes the proof. 
\eproof 

The following corollary of Theorem~\ref{thm:strucsymproj} gives the nearest symmetric matrix of the form $R(X)$ to a given block matrix $A$.

\begin{corollary}
Let $A=\mat{cc} A_1&A_2 \\ A_3 & A_4\rix$, where $A_1,A_2,A_3,A_4 \in \R^{n,n}$. Then  
\begin{equation}
{\rm argmin}_{X \in \C^{n,n},\, X^*=X } {\left \|A-R(X) \right \|}_F =R(\widetilde X),
\end{equation}
where $\widetilde X=\widetilde X_1+i \widetilde X_2$ with $\widetilde X_1=
\frac{{A_1}_H + {A_4}_H}{2}$ and $\widetilde X_2=
\frac{{A_2}_S - {A_3}_S}{2}$.
\end{corollary}

As a corollary of Theorem~\ref{thm:strucsymproj}, we obtain the nearest symmetric PSD matrix of the form $R(X)$ to a given block symmetric matrix $A$. 

\begin{corollary}
Let $A=\mat{cc} A_1&A_2 \\-A_2 & A_1\rix$ be Hermitian, where $A_1,A_2 \in \R^{n,n}$. Then 
\begin{equation}
{\rm argmin}_{X \in \C^{n,n},\, X \succeq 0 } {\left \|A-R(X) \right \|}_F =R(\widetilde X),
\end{equation}
where $\widetilde X= P_{\succeq }(A)$, i.e., PSD projection of $A$.

\end{corollary}

The next corollary of Theorem~\ref{thm:strucsymproj} gives the nearest skew-symmetric matrix of the form $R(X)$ to a given block matrix $A$. 

\begin{corollary}\label{cor:skew-symproj}
Let $A=\mat{cc} A_1&A_2 \\ A_3 & A_4\rix$, where $A_1,A_2,A_3,A_4 \in \R^{n,n}$. Then  
\begin{equation}
{\rm argmin}_{X \in \C^{n,n},\, X^*=-X } {\left \|A-R(X) \right \|}_F =R(\widetilde X),
\end{equation}
where $\widetilde X=\widetilde X_1+i \widetilde X_2$ with $\widetilde X_1=
\frac{{A_1}_S + {A_4}_S}{2}$ and $\widetilde X_2=
\frac{{A_2}_H - {A_3}_H}{2}$.
\end{corollary}


\subsection{ Proof of Theorem~\ref{thm:complex_nhpsdproj}}\label{subsec:proof}

Here, we provide a proof of Theorem~\ref{thm:complex_nhpsdproj} which 
projects a block $2n \times 2n$ matrix onto the set 
$\mathcal{N_S}_\succeq ^{S_{2n}}$  defined in~\eqref{def:struns}. 

\proof First note that for any $X=X_1+iX_2$ where $X_1,X_2 \in \R^{n,n}$, we have 
\begin{eqnarray}
&&{\|A-R(X)\|}_F^2\\
&&={\left\| \mat{cc}A_1 & 
A_2 \\ A_3 & A_4\rix - \mat{cc}X_1 &  X_2 \\ -X_2 & X_1 \rix \right \|}_F^2  \nonumber \\
&=& {\left \| \mat{cc} {A_1}_H +{A_1}_S & \frac{A_2+A_3^T}{2} + \frac{A_2-A_3^T}{2} 
\\  \frac{A_3+A_2^T}{2} + \frac{A_3-A_2^T}{2} &  {A_4}_H +{A_4}_S
 \rix -\mat{cc} {X_1}_H + {X_1}_S & {X_2}_H + {X_2}_S \\
 -{X_2}_H - {X_2}_S & {X_1}_H + {X_1}_S \rix \right  \| }_F^2\nonumber \\ \label{eq:projthm31} \\
 &=& {\left\|  \mat{cc} {A_1}_H &  \frac{A_2+A_3^T}{2} \\ \frac{A_3+A_2^T}{2} &  {A_4}_H  \rix - \mat{cc} {X_1}_H & {X_2}_S \\ -{X_2}_S & {X_1}_H \rix
 \right \|}_F^2  \nonumber \\
 &&+~ {\left\|  \mat{cc} {A_1}_S &  \frac{A_2-A_3^T}{2} \\ \frac{A_3-A_2^T}{2} &  {A_4}_S  \rix - \mat{cc} {X_1}_S & {X_2}_H \\ -{X_2}_H & {X_1}_S \rix
 \right \|}_F^2, \nonumber \\ \label{eq:projthm32}
\end{eqnarray}
where in~\eqref{eq:projthm31} we used the fact that  any $C \in \R^{n,n}$
can be written as $C=C_H+C_S$, and in~\eqref{eq:projthm32} we use that
${\|C\|}_F^2={\|C_H\|}_F^2+{\|C_S\|}_F^2$. Also note that 
for any $X \in \C^{n,n}$ such that $R(X) \in \mathcal{N_S}_\succeq ^{S_{2n}}$
if and only if  $R({X_1}_H + {X_2}_S) $ is PSD. Thus taking the infimum in~\eqref{eq:projthm32} over all $X=X_1 +i X2 \in \C^{n,n}$ such that 
$R(X) \in \mathcal{N_S}_\succeq ^{S_{2n}}$, we obtain that 
\begin{eqnarray}\label{eq:projthm33}
&{\rm argmin}_{X \in \C^{n,n},\, R(X) \in \mathcal{N_S}_\succeq ^{S_{2n}}} 
{\left \| A-R(X)\right \|}_F^2  \nonumber \\
 & ={\rm argmin}_{{X_1}_H, {X_2}_S \in \R^{n,n},\, R({X_1}_H+i {X_2}_S) \succeq 0}   
{\left\|  \mat{cc} {A_1}_H &  \frac{A_2+A_3^T}{2} \\ \frac{A_3+A_2^T}{2} &  {A_4}_H  \rix - \mat{cc} {X_1}_H & {X_2}_S \\ -{X_2}_S & {X_1}_H \rix
 \right \|}_F^2 \nonumber \\
& + {\rm argmin}_{{X_1}_S, {X_2}_H \in \R^{n,n},\, {X_1}_S^T=-{X_1}_S,{X_2}_H^*={X_2}_H } {\left\|  \mat{cc} {A_1}_S &  \frac{A_2-A_3^T}{2} \\ \frac{A_3-A_2^T}{2} &  {A_4}_S  \rix - \mat{cc} {X_1}_S & {X_2}_H \\ -{X_2}_H & {X_1}_S \rix
 \right \|}_F^2. \nonumber \\
\end{eqnarray}
By using Theorem~\ref{thm:strucsymproj}, the first infimum in~\eqref{eq:projthm33} is attained by ${\widetilde X_1}{}_H+ i{\widetilde {X}_2}{}_S$, the PSD projection of $\frac{1}{2}({A_1}_H+{A_4}_H)+i\frac{1}{2}({A_2}_S-{A_3}_S)$. Similarly, from Corollary~\ref{cor:skew-symproj} the second infimum in~\eqref{eq:projthm33} is attained by ${\widetilde X_1}{}_S+ i{\widetilde {X}_2}{}_H$, where ${\widetilde X_1}{}_S=\frac{{A_1}_S+{A_4}_S}{2}$, and 
${\widetilde X_2}{}_H=\frac{{A_2}_H-{A_3}_H}{2}$. This yields from~\eqref{eq:projthm33} that 
\begin{equation}
{\rm argmin}_{X \in \C^{n,n},\, R(X) \in \mathcal{N_S}_\succeq ^{S_{n}} } {\left \|A-R(X) \right \|}_F =
\mat{cc}  {\widetilde X_1}{}_H+ {\widetilde {X}_1}{}_S &
{\widetilde X_2}{}_S+{\widetilde X_2}{}_H \\ -{\widetilde X_2}{}_S-{\widetilde X_2}{}_H  & {\widetilde X_1}{}_H+{\widetilde X_1}{}_S \rix,
\end{equation} 
which concludes the proof. 
\eproof

\newpage

%

\end{document}